\documentclass[12pt,twoside]{amsart}

\date{}
\allowdisplaybreaks[4] \footskip=15pt
\renewcommand{\uppercasenonmath}[1]{}

\numberwithin{equation}{section} \theoremstyle{plain}
\newtheorem*{thm*}{Main Theorem}
\newtheorem{theorem}{Theorem}[section]
\newtheorem{corollary}[theorem]{Corollary}
\newtheorem*{corollary*}{Corollary}
\newtheorem{lemma}[theorem]{Lemma}
\newtheorem*{lemma*}{Lemma}
\newtheorem{proposition}[theorem]{Proposition}
\newtheorem*{proposition*}{Proposition}
\newtheorem{remark}[theorem]{Remark}
\newtheorem*{remark*}{Remark}
\newtheorem{example}[theorem]{Example}
\newtheorem*{example*}{Example}

\newtheorem*{definition*}{Definition}

\newtheorem*{ack*}{ACKNOWLEDGEMENTS}

\newcommand{\pf}{\noindent\begin {proof}}
\newcommand{\epf}{\end{proof}}

\begin{document}
\begin{center}
{\large  \bf Extensions of McCoy Rings}

\vspace{0.8cm} {\small \bf  Zhiling Ying$^{a}$, Jianlong Chen$^{a}$
and
  Zhen Lei$^{a,b}$\\
\vspace{0.6cm} {\rm $^a$ Department of Mathematics, Southeast
University
 \\ Nanjing 210096, P.R. China\\
$^b$Department of Mathematics, Anhui Normal University, Wuhu 241000,
China}\\}
 E-mail: zhilingying@yahoo.com.cn\\
  E-mail:  jlchen@seu.edu.cn\\
  E-mail:  wuhuleizhen@126.com\\

 \end{center}

\bigskip
\centerline { \bf  Abstract}
 \bigskip
\leftskip10truemm \rightskip10truemm

\noindent A ring $R$ is said to be right McCoy if the equation
$f(x)g(x)=0,$ where $f(x)$ and $g(x)$ are nonzero polynomials  of
$R[x],$
 implies that there exists nonzero  $s \in R$
such that $f(x)s  = 0$. It is proven that no proper (triangular)
matrix ring is one-sided McCoy. If  there exists the classical right
quotient ring $Q$ of a ring $R$, then $R$ is right McCoy  if and
only if $Q$ is right McCoy. It is shown that for many polynomial
extensions, a ring $R$ is right McCoy if and only if the polynomial
extension over $R$ is right McCoy. Other basic  extensions of right
McCoy rings are also studied.\leftskip0truemm \rightskip0truemm
\\

\bigskip
\section { \bf Introduction}
\bigskip

\rm Throughout this paper, all rings are associative with identity.
Given a ring $R$,  the polynomial ring over $R$ is denoted by
$R[x]$, and the  ring of $n\times n$ matrices (resp., upper
triangular matrices) over $R$ is denoted by $M_n(R)$ (resp.,
$T_n(R)$).

Recently,  Nielsen \cite{Nielsen06} called  a ring $R$ right McCoy
if the equation $f(x)g(x)=0,$ where $f(x),g(x) \in R[x] \backslash
\{0\},$ implies that there exists $s \in R\backslash \{0\}$ such
that $f(x)s  = 0$. Left McCoy rings are defined analogously. McCoy
rings are the left and right McCoy rings. In \cite[Claim 7 and
8]{Nielsen06}, it is shown that there exists  a left McCoy ring but
not right McCoy. The name ``McCoy" was chosen because McCoy
\cite{McCoy42} had noted that every commutative ring satisfies this
condition. Reversible rings (that is, $ab=0$ implies $ba=0$ for all
$a,b \in R$) are McCoy \cite[Theorem 2]{Nielsen06}, and the
relationships among these rings and other related rings are
discussed in \cite{LCY, Nielsen06}. A ring $R$ is called an
Armendariz ring   \cite{RC97}  if $(\sum_{i=0}^sa_ix^i)(
\sum_{j=0}^tb_jx^j)=0$ in $R[x]$ then $a_ib_j=0$  for all $i$ and
$j$.  Armendariz rings are McCoy by
 definition and an example below shows that
McCoy rings need  not be Armendariz. Therefore, McCoy rings are
shown to be a unifying generalization of reversible rings and
Armendariz rings.

In this paper, at first we consider whether the property ``McCoy" is
Morita invariant. It is proven  that for any ring $R$ and $n\geq 2$,
$M_n(R)$ (resp., $T_n(R)$) is neither left nor right McCoy.
Sequentially,  we argue the property ``McCoy" of some kinds of
polynomial rings.   For many polynomial extensions, a ring $R$ is
right McCoy if and only if the polynomial extension over $R$ is
right McCoy.  It is also proven that if  there exists the classical
right quotient ring $Q$ of a ring $R$, then $R$ is right McCoy  if
and only if $Q$ is right McCoy. Moreover, some examples to answer
questions raised naturally in the process are also given.

\bigskip
\section { \bf Matrix Rings Over  Mccoy Rings}
\bigskip

\ \indent In this section, whether the property ``McCoy" is Morita
invariant and the property ``McCoy" of some subring of upper
triangular matrix ring  are investigated.

We start with  the following fact that is  due to \cite[Theorem
2]{LCY}.
\begin{lemma}   \label{R_n}
A ring $R$ is right (resp., left)  McCoy if and only if the ring
$$R_n=\left\{\left(
\begin{array}{ccccc}
a & a_{12} & a_{13} & \ldots & a_{1n} \\
0 & a      & a_{23} & \ldots & a_{2n} \\
0 & 0      & a      & \ldots & a_{3n} \\
\vdots  & \vdots  & \vdots  & \ddots & \vdots \\
0 & 0      & 0      & \ldots & a
\end{array} \right): a, a_{kl} \in R \right\}$$ is  right (resp., left)  McCoy for any $n \geq 1$.
\end{lemma}

Armendariz rings are McCoy,  but there exists a McCoy  ring which is
 not Armendariz. If $R$ is a McCoy ring,  then  so is $R_4$  by
Lemma \ref{R_n}.  But $R_4$ is not Armendariz by \cite [Example 3]
{KL00}.

Based on Lemma \ref{R_n}, one may suspect that $M_n(R)$ or $T_n(R)$
over a McCoy ring $R$ is still  McCoy. But $M_2(\mathbb{Z}_4)$,
which is not right McCoy ring by \cite{Weiner52}, erases the
possibility. Therefore, the property ``McCoy" is not Morita
invariant. In general, we obtain the following result.
\begin{theorem} \label{matrix(triangular)}
For any   ring  $R$,  $M_n(R)$  (resp., $T_n(R)$) is neither left
nor right  McCoy for any $n> 1$.
\end{theorem}
\begin{proof}
If $n= 2$, then we denote $A=C=e_{12}, B=e_{11},D=-e_{22},$ where
$e_{ij}$'s are the usual matrix units,  and $f(x)=A+Bx, ~ g(x)=C+Dx
\in M_2(R)[x].$ It is clear that  $f(x)g(x)=0$. But if $f(x)P=0$ or
$Qg(x)=0$ for some $P, Q \in M_2(R)$, then $P=Q=0$.

If $n>  2$, then we denote
\begin{eqnarray*}
& F_1(x) &=\left( \begin{array}{cc}
A     & 0 \\
 0      & 0 \\
\end{array} \right) +
\left( \begin{array}{ccc}
B      & 0 \\
 0      & 0 \\
\end{array} \right)x,
%\nonumber
\\
& G_1(x) & =\left( \begin{array}{cc}
C     & 0 \\
 0      & 0 \\
\end{array} \right) +
\left( \begin{array}{ccc}
D      & 0 \\
 0      & I_{n-2} \\
\end{array} \right)x,
%\nonumber
\\
& F_2(x) & =\left( \begin{array}{cc}
A     & 0 \\
 0      & I_{n-2} \\
\end{array} \right) +
\left( \begin{array}{ccc}
B      & 0 \\
 0      & 0 \\
\end{array} \right)x,
%\nonumber
\\
& G_2(x) & =\left( \begin{array}{cc}
C     & 0 \\
 0      & 0 \\
\end{array} \right) +
\left( \begin{array}{ccc}
D      & 0 \\
 0      & 0 \\
\end{array} \right)x  \in M_n(R)[x].
\end{eqnarray*}
So $F_1(x)G_1(x)=0$ and $F_2(x)G_2(x)=0$. But if $SG_1(x)=0$ or
$F_2(x)T=0$ for some $S, T \in M_n(R)$, then $S=T=0$. Therefore,
$M_n(R)$ is neither left nor right McCoy for any $n\geq 2$.

Note that   $f(x), g(x) \in T_2(R)[x],$  and $F_i(x), G_i(x) \in
T_n(R)[x]$ for any $i=1,2$ and $n>  2$.  Similarly,  it is   proven
that $T_n(R)$ is neither left nor right McCoy for any $n\geq 2$.
\end{proof}

\begin{example} \label{eRe example}
A ring  $R$ is  a right (resp., left) McCoy ring if and only if the
ring
$$V(R)=\left\{ \left(
                   \begin{array}{cccccc}
                     a & d & 0 & 0 & 0 & 0 \\
                     0 & b & 0 & 0 & 0 & 0 \\
                     0 & 0 & c & e & 0 & 0 \\
                     0 & 0 & 0 & a & 0 & 0 \\
                     0 & 0 & 0 & 0 & b & f \\
                     0 & 0 & 0 & 0 & 0 & c \\
                   \end{array}
                 \right)
:a,b,c,d,e,f \in R \right\}$$ is  a right (resp., left)  McCoy ring.
\end{example}
\begin{proof}
It is sufficient to prove the case when $R$ is right McCoy. The
other case is similar.

``$\Rightarrow$".  Let $F(x)=\sum_{i=1}^mA_ix^i,
~G(x)=\sum_{j=1}^nB_jx^j$ be nonzero polynomials in $V(R)[x]$ such
that $F(x)G(x)=0$, where
$$A_i=a_{1i}(e_{11}+e_{44})+a_{2i}(e_{22}+e_{55})+a_{3i}(e_{33}+e_{66})+c_{1i}e_{12}+c_{2i}e_{34}+c_{3i}e_{56},$$
$$B_j=b_{1j}(e_{11}+e_{44})+b_{2j}(e_{22}+e_{55})+b_{3j}(e_{33}+e_{66})+d_{1j}e_{12}+d_{2j}e_{34}+d_{3j}e_{56}.$$

Denote $$f_s(x)=\sum_{i=1}^ma_{si}x^i ,~
p_k(x)=\sum_{i=1}^mc_{ki}x^i,$$
$$g_t(x)=\sum_{j=1}^nb_{tj}x^j,~ q_l(x)=\sum_{j=1}^nd_{lj}x^j,$$
where $1\leq s,t,k,l \leq 3$.  Then we have
$$F(x)=f_1(x)(e_{11}+e_{44})+f_2(x)(e_{22}+e_{55})+f_3(x)(e_{33}+e_{66})+p_1(x)e_{12}+p_2(x)e_{34}+p_3(x)e_{56},$$
$$G(x)=g_1(x)(e_{11}+e_{44})+g_2(x)(e_{22}+e_{55})+g_3(x)(e_{33}+e_{66})+q_1(x)e_{12}+q_2(x)e_{34}+q_3(x)e_{56}.$$

To prove that  $V(R)$ is a right McCoy ring, we may choose some
fixed index  $s,t, k$ or $l$ of the set $\{1,2,3\}$  for simplicity
of statement. If $f_s(x)=0$ for some $s,$ then   we can choose
$S=e_{12}$ if $s=1$, $S=e_{56}$ if $s=2$ and $S=e_{34}$ if $s=3$
such that $F(x)S=0.$  Next suppose that $f_s(x)\neq 0$ for any $s$.

Case 1.  $g_t(x)\neq 0$ for some $t$.

 Assume $t=3$. From
 $F(x)G(x)=0$, we have $f_3(x)g_3(x)=0.$ It implies that there exists nonzero
$r_1\in R$ such that $f_3(x)r_1=0$. So $F(x)r_1e_{34}=0$.

Case 2.  $g_t(x)= 0$ for every $t$.

Since $G(x) \neq 0$,  $q_l(x)\neq 0$ for some $l$. Assume
$q_1(x)\neq 0$. Since
 $f_1(x)q_1(x)=0,$ there exists nonzero $r_2 \in R$ such that
 $f_1(x)r_2=0,$   implying  $F(x)r_2e_{12}=0$.

Therefore, $V(R)$  is a right McCoy ring.

``$\Leftarrow$". Assume that  $f(x)g(x)=0,$  where $f(x)=
 \sum_{i=0}^na_ix^{i} \neq 0, g(x) = \sum_{j=0}^mb_jx^{j} \neq 0, a_i,b_j \in
 R.$  Let $F(x)=
 \sum_{i=0}^n(a_iI_6)x^{i}, G(x) = \sum_{j=0}^m(b_jI_6)x^{j} \in
 V(R)[x],$  where $I_6$ is the identity matrix.
Then $F(x)G(x)= [f(x)I_6][g(x)I_6]=0.$ Hence, there exists nonzero
$S=s_1(e_{11}+e_{44})+s_2(e_{22}+e_{55})+s_3(e_{33}+e_{66})+t_1e_{12}+t_2e_{34}+t_3e_{56}
\in V(R)$ such that $F(x)S=0$ because $V(R)$ is McCoy. If $s_i \neq
0$ for some $i \in \{1,2,3\}$ then $f(x)s_i=0$. If $s_i=0$ for every
$i$,  then there exists $t_j \neq0$ for some $j \in \{1,2,3\}$ since
$S\neq0 $.   We also have $f(x)t_j = 0.$  Thus, $R$ is  right McCoy.
\end{proof}

\begin{remark}
(1) For a McCoy ring $R$,  $eRe$ may not be McCoy for some
idempotent $e \in R$.
\\ \indent(2)  $R$ may not be McCoy even if $eRe$ is McCoy for every
nontrivial idempotent $e$ of $R$.
\end{remark}
\begin{proof}
 \indent(1)  Let $R$ be a McCoy ring. Then $V(R)$ in Example \ref{eRe example}  is
 McCoy.  Set $e=e_{11}+e_{22}+e_{44}+e_{55} \in V(R)$. Then  $e$ is an
idempotent of  $V(R)$, but $eV(R)e \cong \big(\begin{smallmatrix}R & R \\
0 & R
\end{smallmatrix}\big)$ is never McCoy  by Theorem
\ref{matrix(triangular)}.

 \indent(2) Let $R= T_2(\mathbb{Z}_2)$.  Clearly,  the nontrivial
idempotents of $R$ are $\big(\begin{smallmatrix}1 & 0 \\ 0 & 0
\end{smallmatrix}\big)$, $\big(\begin{smallmatrix}0 & 0 \\ 0 & 1
\end{smallmatrix}\big)$, $\big(\begin{smallmatrix}1 & 1 \\ 0 & 0
\end{smallmatrix}\big)$ and $\big(\begin{smallmatrix}0 & 1 \\ 0 & 1
\end{smallmatrix}\big)$.  Though  $R$ is not McCoy by Theorem
\ref{matrix(triangular)}, $eRe \cong \mathbb{Z}_2$ is McCoy for
every nontrivial idempotent $e$ of $R$.
\end{proof}

Recall that both reversible rings and Armendariz rings are McCoy and
abelian (i.e., each  idempotent is  central).  So it is natural to
observe the relationships  between them. A ring $R$ is said to be
semi-commutative if $ab=0$ implies $aRb=0$ for $a, b \in R$. Nielsen
showed that semi-commutative (hence abelian)   rings need not be
McCoy in \cite [Section 3] {Nielsen06}.  Conversely, $V(R)$ over
each McCoy ring $R$ is a non-abelian McCoy ring in Example \ref{eRe
example}.

\bigskip
\section { \bf Other extensions of McCoy rings}
\bigskip

\ \indent Basic extensions (including some kinds of polynomial rings
and  classical quotient rings)  of McCoy rings  are investigated in
this section.

\begin{proposition} \label{direct sum}
If  $R_1$ and $R_2$ are right McCoy,  then so is $R=R_1\times R_2$.
\end{proposition}
\begin{proof}
Let  $f(x)=\sum_{i=0}^m(a_i,b_i)x^i$ and
$g(x)=\sum_{j=0}^n(c_j,d_j)x^j \in R[x] \backslash \{0\} $  such
that $f(x)g(x)=0.$ Set $f_1(x)= \sum_{i=0}^ma_ix^i, f_2(x)=
\sum_{i=0}^mb_ix^i, g_1(x)= \sum_{j=0}^nc_jx^j$ and
$g_2(x)=\sum_{j=0}^nd_jx^j.$ Then $f_1(x)g_1(x)=0=f_2(x)g_2(x).$  If
$f_1(x)=0$, then $f(x)(1,0)=0$. If $f_2(x)=0$, then $f(x)(0,1)=0$.
Next suppose $f_1(x) \neq 0$ and $f_2(x) \neq 0$. Since $g(x) \neq
0$, $g_1(x) \neq 0$ or $g_2(x) \neq 0$. If  $g_1(x) \neq 0$, then
there exists nonzero $s_1 \in R_1$ such that $f_1(x)s_1=0$. Thus
$f(x)(s_1, 0)=0$.  If  $g_2(x) \neq 0$, then there exists nonzero
$s_2 \in R_2$ such that $f_2(x)s_2=0$. Thus $f(x)(0, s_2)=0$.
Therefore, $R$ is right McCoy.
\end{proof}

It is natural to ask whether $R$ is a McCoy ring if for any nonzero
proper ideal $I$ of $R$, $R/I$ and $I$ are McCoy, where $I$ is
considered as a McCoy ring without identity. However,  we have a
negative answer to this question by the following example.
\begin{example}
Let $F$ be a field and consider $R=T_2(F),$  which  is not McCoy by
Theorem \ref{matrix(triangular)}. Next we show that $R/I$ and $I$
are McCoy for any  nonzero proper ideal $I$ of $R$. Note that the
only nonzero proper ideals of $R$ are $\big(\begin{smallmatrix}F & F
\\ 0 & 0 \end{smallmatrix}\big)$, $\big(\begin{smallmatrix}0 & F
\\ 0 & F \end{smallmatrix}\big)$ and $\big(\begin{smallmatrix}0 &
F \\ 0 & 0
\end{smallmatrix}\big)$.

First, let $I=\big(\begin{smallmatrix}F & F \\ 0 & 0
\end{smallmatrix}\big)$. Then $R/I \cong F$ and so $R/I$ is McCoy
obviously. Let $f(x)= \sum_{i=0}^m\big(\begin{smallmatrix}a_i & b_i
\\ 0 & 0 \end{smallmatrix}\big)x^i$ and $g(x)= \sum_{j=0}^n\big(\begin{smallmatrix}c_j &
d_j \\ 0 & 0 \end{smallmatrix}\big)x^j$ be nonzero polynomials of
$I[x]$ such that $f(x)g(x)=0,$ implying
$$f_1(x)g_1(x)=f_1(x)g_2(x)=0,            \hspace{4.8cm}       (*)$$ where $f_1(x)=\sum_{i=0}^ma_ix^i,
g_1(x)=\sum_{j=0}^nc_jx^j, g_2(x)=\sum_{j=0}^nd_jx^j \in F[x].$ If
$f_1(x)=0$, then $f(x)e_{11}=0$.  Suppose $f_1(x)\neq 0$. Since
$g(x) \neq 0$, $g_1(x) \neq 0$ or $g_2(x) \neq 0$. From  the
equation $(*)$ and the condition that $F$ is right McCoy,  we have
$f_1(x)s=0$ for some nonzero $s \in F,$ whence $f(x)(se_{11})=0$.
Thus, $I$ is right McCoy,  and $I$ is left  McCoy  since
$e_{12}g(x)=0$.   Next let $J= \big(\begin{smallmatrix}0 & F \\ 0 &
F
\end{smallmatrix}\big)$. Then $R/J$ and $J$ are McCoy by the same
method. Finally, let $K=\big(\begin{smallmatrix}0 & F \\ 0 & 0
\end{smallmatrix}\big)$. Then $R/K\cong F\bigoplus F$ is McCoy.  Since for any $h(x)\in K[x]$, $h(x)e_{12}=e_{12}h(x)=0$, $K$ is
obviously McCoy.
\end{example}

A classical right quotient ring for $R$ is a ring $Q$ which contains
$R$ as a subring in such a way that every regular
  element (i.e., non-zero-divisor) of $R$ is invertible in $Q$ and
$Q=\{ab^{-1}:a, b \in R,~ b ~\mbox{ regular} \}.$  A ring $R$ is
called right Ore if given $a, b \in R$ with $b$ regular there exist
$a_1,  b_1 \in R$ with $b_1$ regular such that $ab_1 = ba_1.$
Classical left quotient rings and left Ore rings are defined
similarly. It is a well-known fact that $R$ is a right (resp., left)
Ore ring if and only if the classical right  (resp., left)  quotient
ring of $R$ exists.

\begin{theorem} \label{classical right quotient  ring}
Suppose that there exists the classical right quotient ring $Q$ of a
ring $R$. Then $R$ is right McCoy if and only if $Q$ is right McCoy.
\end{theorem}
\begin{proof}
``$\Rightarrow$". Let  $F(x)=\sum_{i=0}^m\alpha_ix^i$ and
$G(x)=\sum_{j=0}^n\beta_jx^j$ be nonzero polynomials of $Q[x]$ such
that $F(x)G(x)=0.$  Since $Q$ is a classical right  quotient ring,
we may assume that $\alpha_i = a_iu^{-1},~ \beta_j = b_jv^{-1}$ with
$a_i,~ b_j \in R$ for all $i,j$ and regular elements $u,~v \in R$ by
\cite [Proposition 2.1.16]{McConnell87}. For each $j$, there exist
$c_j \in R$ and a regular element $w \in R$ such that
$u^{-1}b_j=c_jw^{-1}$ also by \cite [Proposition
2.1.16]{McConnell87}. Denote $f_1(x)=\sum_{i=0}^ma_ix^i$ and
$g_1(x)=\sum_{j=0}^nc_jx^j$. Then the equation

$$f_1(x)g_1(x)(vw)^{-1}  = \sum_{i=0}^m \sum_{j=0}^n(a_ic_j)(vw)^{-1}x^{i+j}
=\sum_{i=0}^m \sum_{j=0}^na_i(u^{-1}b_j)v^{-1}x^{i+j}=  F(x)G(x)
=0$$ implies  $f_1(x)g_1(x)=0.$ Thus, there exists a nonzero element
$s \in R$ such that $f_1(x)s=0,$  i.e.,  $a_is=0$ for every $i$.
Then $\alpha_i(us)=a_is=0$ for every $i$. It implies that
$F(x)(us)=0$ and $us$ is a nonzero element of $Q$.  Hence, $Q$ is
right  McCoy.

``$\Leftarrow$". Let  $f(x)=\sum_{i=0}^ma_ix^i$ and
$g(x)=\sum_{j=0}^nb_jx^j \in R[x] \backslash \{0\} $  such that
$f(x)g(x)=0.$   Then  there  exists a nonzero element $\alpha \in Q$
such that $f(x)\alpha=0$ since $Q$ is right McCoy.   Because $Q$ is
a classical right quotient ring,  we can assume $\alpha=au^{-1}$ for
some $a \in R \backslash \{0\}$ and regular element $u$.  Then
$f(x)au^{-1}=f(x)\alpha=0$ implies that $f(x)a=0$.  Therefore, $R$
is a right  McCoy ring.
\end{proof}

By the Goldie Theorem, if $R$ is semiprime left and right Goldie
ring, then $R$ has  the classical left and right  quotient ring.
Hence there exists a class of  rings satisfying the following
hypothesis.

\begin{corollary}
Suppose that there exists the classical left and right quotient ring
$Q$ of a ring $R$. Then $R$ is McCoy  if and only if $Q$ is McCoy.
\end{corollary}

Recall that for a ring $R$ with a ring endomorphism $\alpha:
R\rightarrow R,$  a skew polynomial ring $R[x; \alpha]$  is the ring
obtained by giving the polynomial ring over $R$ with the new
multiplication $xr=\alpha(r)x$ for all $r \in R.$  And a Laurent
polynomial ring  $R[x;x^{-1}]$ is the ring consisting  of all formal
sums $\sum_{i=k}^nr_ix^i$ with obvious addition and multiplication,
where $r_i\in R$ and $k,n$ are (possibly negative) integers. For
rings $R[x]/(x^n)$ and $R[x; \alpha]/(x^n)$, we always consider $n
\geq 2.$

\begin{proposition} \label{skew modulo polynomial}
For a ring $R$ and  an endomorphism $\alpha$ of $R$, the following
statements hold:
\\ \indent (1) $R$ is a right McCoy ring iff $R[x; \alpha]/(x^n)$ is right McCoy.
\\ \indent (2) If $\alpha$ is  monic and  $R$ is a left McCoy ring, then $R[x; \alpha]/(x^n)$ is left McCoy.
\\ \indent (3) If $\alpha^2= \alpha$ and $R$ is a left McCoy ring, then $R[x; \alpha]/(x^n)$ is left McCoy.
\\ \indent (4) If $\alpha$ is an automorphism  and  $R[x; \alpha]/(x^n)$ is left McCoy,  then $R$ is a left McCoy ring.
\end{proposition}
\begin{proof}
Let $F(y)=\sum_{i=0}^pf_iy^i,$ $G(y)=\sum_{j=0}^qg_jy^j $ be nonzero
polynomials in $R[x; \alpha]/(x^n)[y]$  such that $F(y)G(y)=0,$
where $f_i=\sum_{s=0}^{n-1}a_{is}x^s,$
$g_j=\sum_{t=0}^{n-1}b_{jt}x^t \in R[x; \alpha]/(x^n);$ let
$k_s(y)=\sum_{i=0}^pa_{is}y^i$ and $h_t(y)=\sum_{j=0}^qb_{jt}y^j$.
Then
$$[\sum_{s=0}^{n-1}k_s(y)x^s][\sum_{t=0}^{n-1}h_t(y)x^t]=F(y)G(y)=0.   \hspace{5.0cm}(*)$$

(1)  For the ``only if" part, suppose $k_0(y)\neq 0$ and $h_k(y)\neq
0$ with $k$ minimal. Then   $k_0(y)h_k(y)=0$ by the equation $(*)$.
Hence, there exists nonzero $r_1\in R$ such  that $k_0(y)r_1=0,$
implying   $F(y)(r_1x^{n-1})=0.$  If  $k_0(y)=0$, then
$F(y)x^{n-1}=0$.  Therefore, $R[x; \alpha]/(x^n)$ is right McCoy.

For the ``if" part,  let $f(y)=\sum_{i=0}^pa_iy^i,$
$g(y)=\sum_{j=0}^qb_jy^j \in R[y]\backslash \{0\}$ such that
$f(y)g(y)=0.$  Because $f(y)$ and $g(y)$ are nonzero polynomials of
 $R[x; \alpha]/(x^n)[y]$ and  $R[x; \alpha]/(x^n)$ is  right  McCoy, there exists a  nonzero
polynomial $h_1(x)= \sum_{k=0}^{n-1}c_kx^k$  of $R[x; \alpha]/(x^n)$
such that $f(y)h_1(x)=0.$  Let $c_{k_0}\neq 0$ with $k_0$ minimal.
Thus, $f(y)c_{k_0}=0.$ Hence, $R$ is  right McCoy.

(2) If $h_0(y)=0$, then $x^{n-1}G(y)=0$. Next suppose that $h_0(y)
\neq 0$ and $k_l(y)\neq 0$ with $l$ minimal. Thus
$k_l(y)x^lh_0(y)=0, $ implying  that
$k_l(y)[\sum_{j=0}^q\alpha^l(b_{j0})y^j]=0.$  Since $\alpha$ is a
monomorphism,  $\alpha^l(b_{j0})$ are not  all zero, i.e., $
\sum_{j=0}^q\alpha^l(b_{j0})y^j \neq 0$. Because $R$ is left McCoy,
there exists  $r_2 \in R\backslash \{0\}$ such that
$r_2[\sum_{j=0}^q\alpha^l(b_{j0})y^j]=0.$  It implies that
$r_2\alpha^l(b_{j0})=0$  for every $j$, whence
$\alpha^{n-1-l}(r_2)\alpha^{n-1}(b_{j0})=0$ and $\alpha^{n-1-l}(r_2)
\neq 0.$ So $[\alpha^{n-1-l}(r_2)x^{n-1}]G(y)=0.$  It shows that
$R[x; \alpha]/(x^n)$ is left McCoy.

(3) The proof  of (2) needs only  minor modifications to apply here.
If $h_0(y)=0$, then $x^{n-1}G(y)=0$. Next suppose that $h_0(y) \neq
0$ and $k_l(y)\neq 0$ with $l$ minimal. By $\alpha^2= \alpha$ and
$k_l(y)x^lh_0(y)=0, $ we have that
$k_l(y)[\sum_{j=0}^q\alpha(b_{j0})y^j]=0.$ If all
$\alpha(b_{j0})=0,$ then $x^{n-1}G(y)=0.$  If $\alpha(b_{j0})$ are
not  all zero,
 then there exists  $r_3 \in R\backslash \{0\}$ such that
$r_3[\sum_{j=0}^q\alpha(b_{j0})y^j]=0$ since  $R$ is left McCoy. It
implies that $r_3\alpha(b_{j0})=0$  for every $j$, whence
 $(r_3x^{n-1})G(y)=0.$  Thus
$R[x; \alpha]/(x^n)$ is left McCoy.

(4)  Let $f(y)$ and $g(y)$ be the same as the ``if" part in (1).
Using a similar proof, we can obtain that there exist nonzero
$d_{l_0} \in R$ and $0 \leq l_0 \leq n-1$ such that
 $d_{l_0}x^{l_0}g(y)=0,$  i.e., $d_{l_0}\alpha^{l_0}(b_j)=0$ for every
 $j$. Because $\alpha$ is an automorphism,  there exists a nonzero
 element $d'_{l_0} \in R$ such that
 $\alpha^{l_0}(d'_{l_0})=d_{l_0}.$ So $d'_{l_0}b_j=0$ for every
 $j$, whence
 $d'_{l_0}g(y)=0$. Therefore, $R$ is a left McCoy ring.
\end{proof}

\begin{example} \label{endomorphism not auto}
Let $R$ be a left McCoy ring. Consider $R_m (m \geq 2)$  in Lemma
\ref{R_n} and define $\alpha: R_m \rightarrow R_m$ by
$\alpha(A)=aI_m$, where $a$ is the entry on the main diagonal of $A
\in R_m$. Then $R_m$  and $R_m[x; \alpha]/(x^n)$ are left McCoy by
Lemma \ref{R_n}  and Proposition \ref{skew modulo polynomial} (3)
respectively.
\end{example}

From Example \ref{endomorphism not auto}, we know that the condition
``$R[x; \alpha]/(x^n)$ is  left McCoy'' does not imply that $\alpha$
is  monic  or  epic.

\begin{lemma} \label{closed subset}
Let $R$ be a ring and $\Delta$  be a multiplicatively closed subset
of $R$  consisting entirely  of central regular elements. Then $R$
is right  McCoy if and only if $\Delta^{-1}R$ is right McCoy.
\end{lemma}
\begin{proof}
Observe that  it is easy to find a common denominator  for finite
sets of elements in $\Delta^{-1}R$. Using the same way   as  Theorem
\ref{classical right  quotient ring},  the proof is completed.
\end{proof}

\begin{theorem}
For a ring $R$, the following are equivalent:
\\ \indent (1) $R$ is a right McCoy ring.
\\ \indent (2) $R[x]$ is a right McCoy ring.
\\ \indent (3) $R[x;x^{-1}]$ is a right  McCoy ring.
\\ \indent (4) $R[x]/(x^n)$ is a   right McCoy  ring.
\\ \indent (5) $R[\{x_\alpha\}]$ is right  McCoy, where $\{x_\alpha\}$ is
any set of commuting indeterminates over $R$.
\end{theorem}
\begin{proof}
(1) $\Leftrightarrow$ (2)  is due to  \cite [Theorem 1] {LCY} and
 (1) $\Leftrightarrow$ (4) is by Proposition \ref{skew modulo
 polynomial} (1).

(1) $\Rightarrow$  (5). Let $F(y), G(y) \in R[\{x_\alpha\}][y]$ with
$F(y)G(y)=0.$ Then $F(y), G(y) \in R[x_{\alpha_1}, x_{\alpha_2},
\cdots, x_{\alpha_n}][y]$ for some finite subset $\{x_{\alpha_1},
x_{\alpha_2}, \cdots, x_{\alpha_n} \} \subseteq \{x_\alpha\}.$
Following  ``(1) $\Rightarrow$  (2)" and  by induction, the ring
$R[x_{\alpha_1}, x_{\alpha_2}, \cdots, x_{\alpha_n}]$ is right
 McCoy, so there exists nonzero $h_1 \in R[x_{\alpha_1},
x_{\alpha_2}, \cdots, x_{\alpha_n}] \subseteq R[\{x_\alpha\}]$ such
that $F(y)h_1=0$. Hence, $R[\{x_\alpha\}]$ is  right McCoy.

(5) $\Rightarrow$ (1)  is similar to  ``(2) $\Rightarrow$ (1)".

(2) $\Leftrightarrow$ (3). Let $\Delta = \{1,x, x^2,\cdots \}$. Then
clearly $\Delta$ is a multiplicatively closed subset of $R[x]$
consisting entirely  of central regular elements.  Since $R[x;
x^{-1}] = \Delta^{-1}R[x],$  $R[x; x^{-1}]$ is right McCoy iff
$R[x]$ is right McCoy by Lemma \ref{closed subset}.
\end{proof}

According to \cite{Krempa96},  an endomorphism $\alpha$ of a ring
$R$ is said to be rigid if $a\alpha(a) = 0$ implies $a = 0$ for
every $a \in R$. Later, Hong  et al.  called  a ring $R$  an
$\alpha$-rigid ring \cite{Hong00}  if there exists a rigid
endomorphism $\alpha$ of $R.$ Clearly, if $R$ is an  $\alpha$-rigid
ring, then $\alpha$ is a monomorphism
 and $R$ is reduced (hence  McCoy). Combining
 Proposition \ref{skew modulo polynomial} (1) and  (2), we obtain
 \begin{corollary}
If  $R$ is  an $\alpha$-rigid  ring,  then  $R[x; \alpha]/(x^n)$ is
a  McCoy ring.
\end{corollary}

If $R$ is an $\alpha$-rigid ring, then $R[x; \alpha]$ is McCoy
because $R[x; \alpha]$ is a reduced ring by \cite[Corollary
3.4]{Krempa96} or \cite [Proposition 5]{Hong00}.  In general, $R[x;
\alpha]$ may not be  McCoy even if $R$ is a commutative reduced ring
and $\alpha$ is an automorphism of $R$.  To show it, we use a ring
given in \cite [Example 6]{KL00}.

\begin{example} \label{skew polynomial}
Let $R=\mathbb{Z}_2 \bigoplus \mathbb{Z}_2$ and $\alpha:
R\rightarrow R$  defined by $\alpha((a,b))=(b,a)$. Then both $R$ and
$R[x; \alpha]/(x^n)$  are McCoy, but $R[x; \alpha]$ is neither left
nor right McCoy.
\end{example}
\begin{proof}
$R$ is a McCoy ring since $R$ is commutative, and so is $R[x;
\alpha]/(x^n)$  by Proposition \ref{skew modulo polynomial} (2)
since $\alpha$ is an automorphism of $R$.  Let
$f(y)=(1,0)+[(1,0)x]y$ and $g(y)=(0,1)+[(1,0)x]y$ be elements in
$R[x; \alpha][y].$ Then $f(y)g(y)=0.$ Denote
$h_1(x)=\sum_{i=0}^m(a_i,b_i)x^i, h_2(x)=\sum_{j=0}^n(c_j,d_j)x^j
\in R[x; \alpha].$   If $f(y)h_1(x)=0,$ then
$\sum_{i=0}^m(a_i,0)x^i+[\sum_{i=0}^m(b_i,0)x^{i+1}]y=0,$ whence
$a_i=b_i=0$ for all $i$.  Thus $h_1(x)=0.$  If $h_2(x)g(y)=0,$ then
$\sum_{j=0}^n(c_j,d_j)\alpha^j((0,1))x^j+[\sum_{j=0}^n(c_j,d_j)\alpha^{j}((1,0))x^{j+1}]y=0,$
whence $c_j=d_j=0$ for all  $j$. Thus $h_2(x)=0.$ Therefore, $R[x;
\alpha]$ is neither left nor right McCoy.
\end{proof}

 \centerline {\bf Acknowledgments}
 This research was supported  by
the National Natural Science Foundation of China (10571026), the
Natural Science Foundation of Jiangsu Province (2005207),  and the
Specialized Research Fund for the Doctoral Program of Higher
Education (20060286006).

\end{document}